\documentclass[fleqn]{mat01}
\usepackage{times,mathtimy,amssymb,latexsym}
\begin{document}

\setcounter{page}{49} \firstpage{49}

\newtheorem{theore}{Theorem}
\renewcommand\thetheore{\arabic{theore}}
\newtheorem{theor}[theore]{\bf Theorem}
\newtheorem{defini}{\rm DEFINITION}
\newtheorem{lem}{Lemma}
\newtheorem{case}{\it Case}

\renewcommand\theequation{\arabic{section}.\arabic{equation}}

\title{On the problem of isometry of a hypersurface preserving
mean curvature}

\markboth{H\"{u}lya Ba\u{g}datli and Ziya Soyu\c{c}ok}{Isometry of
a hypersurface preserving mean curvature}

\author{H\"{U}LYA BA\u{G}DATLI$^{1}$ and ZIYA SOYU\c{C}OK$^{2}$}

\address{$^1$Department of Mathematics, Marmara University,
Istanbul, Turkey\\
\noindent $^2$Department of
Mathematics, Y\i ld\i z Technical University, Istanbul, Turkey\\
\noindent E-mail: hbagdatli@marmara.edu.tr;
zsoyucok@yildiz.edu.tr}

\volume{117}

\mon{February}

\parts{1}

\pubyear{2007}

\Date{MS received 26 October 2004}

\begin{abstract}
The problem of determining the {\it Bonnet hypersurfaces in}
$R^{n+1}$, for $n>1$, is studied here. These hypersurfaces are by
definition those that can be isometrically mapped to another
hypersurface or to itself (as locus) by at least one nontrivial
isometry preserving the mean curvature. The other hypersurface
and/or (the locus of) itself is called {\it Bonnet associate} of
the initial hypersurface.

The orthogonal net which is called \hbox{\it $A$-net} is special
and very important for our study and it is described on a
hypersurface. It is proved that, non-minimal hypersurface in
$R^{n+1}$ with no umbilical points is a Bonnet hypersurface if and
only if it has an\break $A$-net.
\end{abstract}


\maketitle

\section{Introduction}

The isometry problem or more specifically the isometrical
deformation problem of surfaces in a 3-dimensional Euclid space
form preserving mean curvature has been studied by a number of
mathematicians. One of the first mathematicians who has studied
this subject is Bonnet \cite{2}. Bonnet, after pointing out that
the isometry problem of a surface preserving mean curvature can
not be solved in the general case, showed that all surfaces with
constant mean curvature can be isometrically mapped to each other
and deformable surfaces with non constant mean curvature are
isothermic Weingarten surfaces which can be deformable to the
revolution surfaces. Cartan \cite{3} considered the problem as an
application of differential forms and after rather long
calculations, classified Weingarten surfaces which can be
deformable to the revolution surfaces and showed that there are
finite number surfaces in each class.

Recently, this problem has been reconsidered. Chern \cite{4},
using differential forms, obtained a characterization for an
isometrical deformation preserving mean curvature. Roussos
\cite{11} who has many works on this subject, using Chern's
method, in the general case, obtained a characterization for the
isometry preserving mean curvature. Some of the important
contributions to this subject are by Voss \cite{15}, Bobenko and
Either \cite{1}, Kenmotsu \cite{8}, Colares and Kenmotsu \cite{6},
Roussos \cite{12}, Xiuxiong and Chia-Kuei\break \cite{16}.

The surface which admits an isometry-preserving mean curvature is
called Bonnet surface. One of the Bonnet surfaces which can be
isometrically mapped to each other is called Bonnet associate of
the other surface. Similar definitions are valid for the Bonnet
hypersurfaces.

Kokubu \cite{10} considered deformable Bonnet hypersurfaces and
showed that they consist of Bonnet surfaces in $R^{3}$ and Bonnet
hypersurfaces in $R^{4}$.

Later, Soyu\c{c}ok \cite{13} examined the problem of determining
Bonnet surfaces in $R^{3}$ and showed that a surface is a Bonnet
surface if and only if a surface has a special isothermal
parameter system. Moreover, as an application of this, he obtained
Cartan's results and explained geometrically why there are a
finite number of surfaces in the classes determined by Cartan.

Soyu\c{c}ok \cite{14}, in another work, examined the problem
determining Bonnet hypersurfaces in $R^{4}$ and obtained
quantities of these hypersurfaces. Furthermore, he showed that a
3-dimensional hypersurface is a Bonnet hypersurface if and only if
it has a special orthogonal net.

\section{Preliminaries}

\setcounter{equation}{0}

Let us consider $M$ and $M^{\prime}$ to be non-minimal Bonnet
associate hypersurfaces with no umbilical points in $R^{n+1}$.
Therefore, the ranks of their shape operators are the same and
less than 3 \cite{10}.

Let the local orthonormal frame fields of $M$ and $M^{\prime}$
consist of the principal vectors $e_{1},e_{2},\dots,e_{n}$ and
$e_{1}^{\prime},e_{2}^{\prime},\dots,e_{n}^{\prime}$ where
corresponding principal curvatures are $k_{1},k_{2},(k_{1}\neq
k_{2}),k_{3}(=0),k_{4}(=0),\dots,k_{n}(=0)$ and
$k_{1}^{\prime},k_{2}^{\prime},(k_{1}^{\prime}\neq
k_{2}^{\prime}),k_{3}^{\prime}(=0),k_{4}^{\prime}(=0),\dots,
k_{n}^{\prime}(=0)$ respectively. Since, the null spaces of $M$
and $M^{\prime}$ are the same at the corresponding points
\cite{10},
\begin{align*}
&\text{span}\{e_{1},e_{2}\} =
\text{span}\{e_{1}^{\prime},e_{2}^{\prime}\},\\[.2pc]
&\text{span}\{e_{3},e_{4},\dots,e_{n}\} =
\text{span}\{e_{3}^{\prime},e_{4}^{\prime},\dots,
e_{n}^{\prime}\}.
\end{align*}

Because $M$ and $M^{\prime}$ are associate to each other, the
first and the second of the curvature lines of $M^{\prime}$
correspond to the orthogonal curves on $M$. These curves are
called Bonnet curves of $M$.

Let $E_{1}$ and $E_{2}$ be the unit vectors of the bisectors of
Bonnet curves with curvature lines (first and second) of $M$. We
take $\{E_{1},E_{2},E_{3}=e_{3},\dots,E_{n}=e_{n}\}$ as the
orthonormal frame field for $M$ and we can also choose
$\{E_{1},E_{2},E_{3},\dots,E_{n}\}$ as the orthonormal frame field
for $M^{\prime}$.

Let us denote the shape operator and the Riemannian connection of
$M$ by $L$ and $D$ respectively. Then
\begin{align}
\hskip -4pc L=\frac{1}{2}\left(\begin{array}{ccccc}
k_{1}(1+\cos\theta)+k_{2}(1-\cos \theta) &(k_{1}-k_{2})\sin \theta
&0 &\dots &0\\[.2pc]
(k_{1}-k_{2})\sin \theta &k_{1}(1-\cos\theta)+k_{2}(1+\cos \theta)
&0 &\dots &0\\[.2pc]
0 &0 &0 &\dots &0\\[.2pc]
\dots &\dots &\dots &\dots &\dots\\[.2pc]
0 &0 &0 &\dots &0
\end{array}
\right)\nonumber\\
\end{align}
can be written to the frame field $\{E_{1},E_{2},E_{3},
\dots,E_{n}\}$ where $\theta$ is the angle between $e_{1}$ and
$e_{1}^{\prime}$.

Because the angle between $e_{1}^{\prime}$ and $E_{1}$ is
$-\frac{\theta }{2}$, the shape operator of $M^{\prime}$ with
respect to the same frame field can be written as
\begin{align}
\hskip -4pc L^{\prime}=\frac{1}{2}\left(\begin{array}{ccccc}
k_{1}(1\!+\!\cos\theta)+k_{2}(1\!-\!\cos\theta)
&-(k_{1}-k_{2})\sin \theta
&0 &\dots &0\\[.2pc]
-(k_{1}-k_{2})\sin \theta
&k_{1}(1\!-\!\cos\theta)+k_{2}(1\!+\!\cos\theta)
&0 &\dots &0\\[.2pc]
0 &0 &0 &\dots &0\\[.2pc]
\dots &\dots &\dots &\dots &\dots\\[.2pc]
0 &0 &0 &\dots &0
\end{array}
\right).\nonumber\\
\end{align}

We shall use the following convention on the ranges of indices in
this work:
\begin{align*}
1 &\leq i,j,k,\ldots \leq 2,\\[.2pc]
3 &\leq p,q,r,\ldots \leq n,\\[.2pc]
1 &\leq A,B,C,\ldots \leq n.
\end{align*}

\section{Fundamental equations}

\setcounter{equation}{0}

\subsection{\it Codazzi equations}

Let us define the functions $k,\bar{k}$ and $t$ as follows:
\begin{align}
k &= \frac{1}{2}\{k_{1}(1+\cos \theta)+k_{2}(1-\cos \theta)\}
=H+J\cos \theta,\nonumber\\[.4pc]
\bar{k} &= \{k_{1}(1-\cos \theta)+k_{2}(1+\cos
\theta)\}=H-J\cos \theta,\nonumber\\[.4pc]
t &= \frac{1}{2}(k_{1}-k_{2})\sin \theta =J\sin \theta,
\end{align}
where $H$ is a mean curvature and
\begin{equation}
H=\frac{k_{1}+k_{2}}{2}\neq 0,\quad J=\frac{k_{1}-k_{2}}{2}\neq 0.
\end{equation}

The Codazzi equation \cite{9} is given by
\begin{equation}
(D_{E_{i}}^{L}) E_{j}=(D_{E_{j}}^{L}) E_{i}.
\end{equation}

Accordingly, using the shape operator of $M$ and the above
definitions we can write the Codazzi equations of $M$ in the form:
\begin{align}
&E_{2}(k)+(k-\bar{k})w_{2}^{1}(E_{1})=E_{1}(t)+2tw_{1}^{2}(E_{2}),\nonumber\\[.2pc]
&E_{1}(\bar{k})+(\bar{k}
-k)w_{1}^{2}(E_{2})=E_{2}(t)+2tw_{2}^{1}(E_{1}),\nonumber\\[.2pc]
&\bar{k}w_{2}^{p}(E_{1})-kw_{1}^{p}(E_{2})=t[w_{2}^{p}(E_{2})-w_{1}^{p}(E_{1})],\nonumber\\[.2pc]
&E_{p}(k)+kw_{p}^{1}(E_{1})=t[2w_{1}^{2}(E_{p})-w_{p}^{2}(E_{1})],\nonumber\\[.2pc]
&E_{p}(t)+tw_{p}^{1}(E_{1})=\bar{k}[w_{1}^{2}(E_{p})-w_{p}^{2}(E_{1})] -kw_{1}^{2}(E_{p}),\nonumber
\end{align}
\begin{align}
&E_{p}(t)+tw_{p}^{2}(E_{2})=k[w_{2}^{1}(E_{p})-w_{p}^{1}(E_{2})]
-\bar{k}w_{2}^{1}(E_{p}),\nonumber\\[.2pc]
&E_{p}(\bar{k})+\bar{k}w_{p}^{2}(E_{2})=t[2w_{2}^{1}(E_{p})-w_{p}^{1}(E_{2})],\nonumber\\[.2pc]
&kw_{1}^{p}(E_{p})=-tw_{2}^{p}(E_{p}),\nonumber\\[.2pc]
&kw_{1}^{q}(E_{p})=-tw_{2}^{q}(E_{p}),\nonumber\\[.2pc]
&\bar{k}w_{2}^{p}(E_{p})=-tw_{1}^{p}(E_{p}),\nonumber\\[.2pc]
&k[w_{p}^{1}(E_{q})-w_{q}^{1}(E_{p})] =-t[w_{p}^{2}(E_{q})-w_{q}^{2}(E_{p})],\nonumber\\[.2pc]
&\bar{k}[w_{p}^{2}(E_{q})-w_{q}^{2}(E_{p})]
=-t[w_{p}^{1}(E_{q})-w_{q}^{1}(E_{p})],
\end{align}
where $w_{B}^{A}$ are the connection forms.

Because of (2.1), (2.2) and (3.1), in order to find the Codazzi
equations of $M^{\prime}$ it is enough to replace $t$ by $-t$ in
(3.4). Therefore, from the Codazzi equations of $M$ and
$M^{\prime}$ the following equations are obtained:
\begin{align}
&w_{1}^{p}(E_{p})=w_{2}^{p}(E_{p})=0,\nonumber\\[.2pc]
&w_{1}^{q}(E_{p})=w_{2}^{q}(E_{p})=0,\nonumber\\[.2pc]
&w_{p}^{2}(E_{2})=w_{p}^{1}(E_{1}),\nonumber\\[.2pc]
&w_{p}^{1}(E_{2})=-w_{p}^{2}(E_{1})=2w_{2}^{1}(E_{p}),\nonumber\\[.2pc]
&Hw_{p}^{1}(E_{2})=0,\nonumber\\[.2pc]
&E_{p}(k)+kT_{p}=0,\nonumber\\[.2pc]
&E_{p}(\bar{k})+\bar{k}T_{p}=0,\nonumber\\[.2pc]
&E_{p}(t)+tT_{p}=0,\nonumber\\[.2pc]
&E_{2}(k)+(k-\bar{k})h=0,\nonumber\\[.2pc]
&E_{1}(\bar{k})+(\bar{k}-k)\bar{h}=0,\nonumber\\[.2pc]
&E_{1}(t)+2t\bar{h}=0,\nonumber\\[.2pc]
&E_{2}(t)+2th=0,
\end{align}
where the functions $T_{p}$, $h$ and $\bar{h}$ are defined as
\begin{align}
&w_{p}^{2}(E_{2})=w_{p}^{1}(E_{1})=T_{p},\nonumber\\[.2pc]
&h=w_{2}^{1}(E_{1}), \quad \bar{h}=w_{1}^{2}(E_{2}).
\end{align}

Since we assume that our hypersurface is a non-minimal one, $H\neq
0$ and from the fifth equation of (3.5) we get
\begin{equation*}
w_{p}^{1}(E_{2})=0.
\end{equation*}
From the fourth equation of (3.5),
\begin{equation}
w_{p}^{1}(E_{2})=w_{p}^{2}(E_{1})=w_{2}^{1}(E_{p})=0
\end{equation}
are obtained.\newpage

Therefore, we have the following lemma.

\begin{lem}
A Bonnet hypersurface must satisfy the Codazzi equations{\rm ,}
\begin{align}
&E_{p}(k)+kT_{p}=0,\nonumber\\[.2pc]
&E_{p}(\bar{k})+\bar{k}T_{p}=0,\nonumber\\[.2pc]
&E_{p}(t)+tT_{p}=0,\nonumber\\[.2pc]
&E_{2}(k)+(k-\bar{k})h=0,\nonumber\\[.2pc]
&E_{1}(\bar{k})+(\bar{k}-k)\bar{h}=0,\nonumber\\[.2pc]
&E_{1}(t)+2t\bar{h}=0,\nonumber\\[.2pc]
&E_{2}(t)+2th=0.
\end{align}
\end{lem}

\subsection{\it Gauss equations}

The Gauss equations of $M$ hypersurface \cite{9} are given by the
equations
\begin{align*}
&D_{E_{A}} (D_{E_{B}}^{E_{C}}) -D_{E_{B}} (D_{E_{A}}^{E_{C}})
-D_{[E_{A},E_{B}]}^{E_{C}}\\[.3pc]
&\quad\, = \langle L(E_{B}),E_{C}\rangle L(E_{A})-\langle
L(E_{A}),E_{C}\rangle L(E_{B}).
\end{align*}

Using (3.7) and (3.8) in the Gauss equations of $M$, we have the
following lemma.

\begin{lem}
The Gauss equations of $M$ can be written as follows{\rm :}
\begin{align}
&E_{1}(\bar{h})+E_{2}(h)+\sum_{p=3}^{n}T_{p}^{2}+h^{2}+ \bar{h}^{2}=t^{2}-k\bar{k},\nonumber\\[.2pc]
&E_{p}(h)+T_{p}h=0,\nonumber\\[.2pc]
&E_{p}(\bar{h})+T_{p}\bar{h}=0,\nonumber\\[.2pc]
&E_{p}(T_{q})-w_{q}^{r}(E_{p})T_{r}+T_{p}T_{q}=0,\nonumber\\[.2pc]
&E_{1}(T_{p})-w_{p}^{q}(E_{1})T_{q}=0,\nonumber\\[.2pc]
&E_{2}(T_{p})-w_{p}^{q}(E_{2})T_{q}=0,\nonumber\\[.2pc]
&E_{p}(w_{r}^{s}(E_{q}))
+w_{r}^{t}(E_{q})w_{t}^{s}(E_{p})-E_{q}(w_{r}^{s}(E_{p}))\nonumber\\[.2pc]
&\quad\, -w_{r}^{t}(E_{p})w_{t}^{s}(E_{q})-w_{\!p}^{t}(E_{q})
w_{r}^{s}(E_{t})+w_{q}^{t}(E_{p})w_{r}^{s}(E_{t})=0.
\end{align}
\end{lem}

\subsection{\it Structure equations}

Let $\{w^{1},w^{2},\dots,w^{n}\}$ be the dual frame field of
$\{E_{1},E_{2},\dots,E_{n}\}$. Thus, by means of (3.6) and (3.7)
\begin{align*}
w_{2}^{1} &= hw^{1}-\bar{h}w^{2},\\[.2pc]
w_{p}^{1} &= T_{p}w^{1},\\[.2pc]
w_{p}^{2} &= T_{p}w^{2},\\[.2pc]
w_{p}^{q} &= w_{p}^{q}(E_{A})w^{A}
\end{align*}
are obtained. Accordingly, the structure equations $\hbox{d}w^{L}=
-w_{J}^{L}\wedge w^{J}$ can be written as follows:
\begin{align}
\hbox{d}w^{1} &= -hw^{1}\wedge w^{2}-T_{p}w^{1}\wedge w^{p},\nonumber\\[.2pc]
\hbox{d}w^{2} &= \bar{h}w^{1}\wedge w^{2}-T_{p}w^{2}\wedge
w^{p},\nonumber\\[.2pc]
\hbox{d}w^{q} &= -w_{p}^{q}(E_{A})w^{A}\wedge w^{p}.
\end{align}

Now, we consider the differential $\hbox{d}f=f_{A}w^{A}$ of $f$
with respect to the dual frame field
$\{w^{1},w^{2},\dots,w^{n}\}$. Using (3.7) and (3.10) in
$\hbox{d}^{2}f=0$, we obtain the compatibility equations
\begin{align}
&f_{12}+hf_{1}-f_{21}-\bar{h}f_{2}=0,\nonumber\\[.2pc]
&f_{1p}+T_{p}f_{1}-f_{p1}+f_{q}w_{p}^{q}(E_{1})=0,\nonumber\\[.2pc]
&f_{2p}+T_{p}f_{2}-f_{p2}+f_{q}w_{p}^{q}(E_{2})=0,\nonumber\\[.2pc]
&f_{rp}+T_{p}f_{r}-f_{pr}+f_{q}w_{p}^{q}(E_{r})=0.
\end{align}

The sixth and seventh equations of the system (3.7) can be written
in the form
\begin{equation}
(\ln \vert t\vert)_{1}=-2\bar{h},\quad (\ln \vert t\vert)_{2}=-2h.
\end{equation}

So, using the first equation of (3.11),
\begin{equation}
h_{1}=\bar{h}_{2}.
\end{equation}

\section{New coordinates}

\setcounter{equation}{0}

Let $g_{AB}$ be the components of the first fundamental tensor.
The form $w^{A}$ can be written as
\begin{equation*}
w^{A}=\sqrt{g_{AA}}\hbox{d}x^{A}\quad (\hbox{not sum in} \ A)
\end{equation*}
in the suitable system of local coordinates
$(x^{1},x^{2},\dots,x^{n})$ \cite{5}. So, from the structure
equations (3.10), we get
\begin{align}
&w_{p}^{q}(E_{1})=0, \quad w_{p}^{q}(E_{2})=0,\\[.2pc]
&w_{r}^{A}(E_{p})=0,\\[.2pc]
&g_{pp}=c^{p}(x^{p}),\\[.2pc]
&h=\frac{(\sqrt{g_{11}}) _{x^{2}}}{\sqrt{g_{11}g_{22}}}=(\ln
\sqrt{g_{11}})_{2},\quad \bar{h}=
\frac{(\sqrt{g_{22}})_{x^{1}}}{\sqrt{g_{11}g_{22}}}=(\ln
\sqrt{g_{22}})_{1},\\[.2pc]
&T_{p}=\frac{(\sqrt{g_{11}})_{x^{p}}}{\sqrt{g_{11}g_{pp}}}=
\frac{(\sqrt{g_{22}})_{x^{p}}}{\sqrt{g_{22}g_{pp}}},
\end{align}
where
\begin{equation}
f_{x^{A}}=\frac{\partial f}{\partial x^{A}}.
\end{equation}

Because of (4.4) in (3.12), equation (3.13) is equivalent to
\begin{equation*}
\left(\ln \frac{g_{11}}{g_{22}}\right)_{x^{1}x^{2}}=0.
\end{equation*}

Therefore, we have
\begin{equation}
\frac{g_{11}}{a(x^{1})}=\frac{g_{22}}{b(x^{2})},
\end{equation}
where $a(x^{1})$ and $b(x^{2})$ are arbitrary functions. Let us
introduce the new coordinates by means of the following scaling
transformation:
\begin{equation*}
\bar{x}^{1}=\int \sqrt{a(x^{1})}\hbox{d}x^{1},
\quad\bar{x}^{2}=\int \sqrt{b(x^{2})}\hbox{d}x^{2}, \quad
\bar{x}^{p}=\int \sqrt{c^{p}(x^{p})}\hbox{d}x^{p}.
\end{equation*}

Then, from (4.7) and (4.3) we get
\begin{equation}
g_{11}=g_{22},\quad g_{pp}=1.
\end{equation}

Here we have again denoted the new coordinates by $x^{1},
x^{2},\dots,x^{n}$. Accordingly (4.5) reduces to
\begin{equation}
T_{p}=\frac{(\sqrt{g_{11}})_{x^{p}}}{\sqrt{g_{11}}}=
\frac{(\sqrt{g_{22}}) _{x^{p}}}{\sqrt{g_{22}}}.
\end{equation}

Because of (4.1) the fifth and sixth equations of the Gauss
equations (3.9) take the forms
\begin{equation*}
E_{1}(T_{p})=0 \quad \hbox{and} \quad E_{2}(T_{p})=0.
\end{equation*}

So, $T_{p}$ does not depend on $x^{1}$ and $x^{2}$.

Moreover, using (4.2), the fourth equation of (3.8) can be written
as
\begin{equation}
E_{p}(T_{q})+T_{p}T_{q}=0
\end{equation}
and so
\begin{equation}
E_{p}(T_{p})+(T_{p})^{2}=0.
\end{equation}

Now, let us consider the following two cases:

\begin{enumerate}
\renewcommand\labelenumi{{\it Case}~\arabic{enumi}.}
\leftskip 1.65pc
\item Some $T_{p}$'s are not zero.

\item All $T_{p}$'s are zero.\vspace{-1pc}
\end{enumerate}

\begin{case}{\rm Assume that some $T_{p}$'s are not zero. Let the
number of non zero $T_{p}$'s be $m-2$. Renewing the indices $p$ of
non zero $T_{p}$'s we can choose them as $3,4,\dots,r$. Hence we
can take the non zero $T_{p}$'s to be $T_{3},T_{4},\dots,T_{r}$
$(3\leq r\leq n)$.

Solving the equations (4.10) and (4.11), we find
\begin{equation}
T_{p}=\frac{C^{p}}{C^{3}x^{3}+\dots +C^{p}x^{p}+\dots
+C^{r}x^{r}},
\end{equation}
where $C^{3},\dots,C^{r}$ are non zero constants of integration.

Accordingly, from (4.9) we have
\begin{equation}
g_{11}=g_{22}=(C^{3}x^{3}+C^{4}x^{4}+\dots +C^{r}x^{r})^{2}\xi
(x^{1},x^{2}),
\end{equation}
where $\xi (x^{1},x^{2})$ is an arbitrary function.

By means of (4.4) and (4.13), the second and the third Gauss
equations of (3.9) are automatically satisfied

Using (4.8) and (4.9), let us solve the third\ equation of (3.7)
and then let us solve eqs.~(3.12), (4.4) and (4.13). Thus
\begin{equation}
t=\frac{C}{(C^{3}x^{3}+\dots +C^{r}x^{r})\xi (x^{1},x^{2})},
\end{equation}
where $C$ is a constant.

On the other hand, using (4.8) and (4.12), solutions of the first
and the second equations of the system (3.7) are obtained as
\begin{equation}
k=\frac{\Psi (x^{1},x^{2})}{(C^{3}x^{3}+\dots +C^{r}x^{r})},\quad
\bar{k}=\frac{\bar{\Psi} (x^{1},x^{2})}{(C^{3}x^{3}+\dots
+C^{r}x^{r})}.
\end{equation}

Now, let us find the mean curvature $H$. Since
$H=\frac{k+\bar{k}}{2}$, according to (3.1) we have
\begin{equation}
H=\frac{\mathfrak{H}(x^{1},x^{2})}{(C^{3}x^{3}+\dots
+C^{r}x^{r})}.
\end{equation}

Moreover, from (3.1), (4.15) and (4.16), we have
\begin{equation}
\theta =\theta (x^{1},x^{2})
\end{equation}
and
\begin{equation}
J=\frac{\mathfrak{J}(x^{1},x^{2})}{(C^{3}x^{3}+\dots
+C^{r}x^{r})}.
\end{equation}

All of Codazzi's equations~(3.8) are satisfied except the fourth
and the fifth equation and all of Gauss' equations~(3.8) are
satisfied except the first equation. Now, we have to consider
these three equations.

By using (3.1), (4.4), (4.6), (4.14), (4.16), (4.17) and (4.18),
the fourth and the fifth equations of Codazzi's equations~(3.8)
can be written as follows:
\begin{equation}
\frac{\mathfrak{H}_{x^{2}}}{\mathfrak{J}}=\frac{\theta_{x^{2}}}{\sin
\theta}, \quad \frac{\mathfrak{H}_{x^{1}}}
{\mathfrak{J}}=-\frac{\theta_{x^{1}}}{\sin \theta}.
\end{equation}

The first equation of Gauss's equations (3.8) is reduced to the
following equation:
\begin{align}
&\vert \mathfrak{J}\vert \sin \theta \{(\ln \vert
\mathfrak{J}\vert \sin \theta)_{x^{1}x^{1}}+ (\ln \vert
\mathfrak{J}\vert \sin \theta)_{x^{2}x^{2}}\}
-2\{(C^{3})^{2}+\dots + (C^{r})^{2}\}\nonumber\\[.2pc]
&\quad\, = 2 (\mathfrak{H}^{2}-\mathfrak{J}^{2}).
\end{align}

\section{The fundamental theorem}

\setcounter{equation}{0}

Let us denote the components of the second fundamental tensor with
$b_{AB}$. They are given by \cite{7}
\begin{equation*}
b_{AB}=\left\langle L\left(\frac{\partial }{\partial
x^{A}}\right),\frac{\partial }{\partial x^{B}}\right\rangle.
\end{equation*}
Thus, we get
\begin{equation*}
b_{11}=g_{11}k,\quad b_{12}=\sqrt{g_{11}g_{22}}t, \quad
b_{22}=g_{22}\bar{k},\quad b_{iq}=b_{pq}=0.
\end{equation*}

By using (4.13), (4.14) and (4.15), non zero components can be
written as
\begin{align*}
b_{11} &= (C^{3}x^{3}+\dots +C^{r}x^{r})
\frac{(\mathfrak{H}+\mathfrak{J} \sin \theta)}{\vert
\mathfrak{J}\vert \sin \theta},\\[.3pc]
b_{22} &= (C^{3}x^{3}+\dots +C^{r}x^{r})
\frac{(\mathfrak{H}-\mathfrak{J} \sin \theta)}{\vert
\mathfrak{J}\vert \sin \theta},\\[.3pc]
b_{12} &= C(C^{3}x^{3}+\dots +C^{r}x^{r}).
\end{align*}

We can make a suitable scale transformation $C=\in
=sgn\mathfrak{J}$.

Thus, the components of the first and the second fundamental
tensors of Bonnet hypersurface $M$ in an orthogonal coordinate
system $x^{1},x^{2},\dots,x^{n}$ are obtained as
\begin{align}
g_{11} &= g_{22}=(C^{3}x^{3}+\dots +C^{r}x^{r})^{2} \frac{1}{\vert
\mathfrak{J}\vert\sin \theta}, \quad g_{pp}=1\nonumber\\[.3pc]
b_{11} &= (C^{3}x^{3}+\dots +C^{r}x^{r})
\frac{(\mathfrak{H}+\mathfrak{J} \sin \theta)}{\vert
\mathfrak{J}\vert \sin \theta},\nonumber\\[.3pc]
b_{22} &= (C^{3}x^{3}+\dots +C^{r}x^{r})
\frac{(\mathfrak{H}-\mathfrak{J} \sin \theta)}{\vert
\mathfrak{J}\vert \sin \theta},\nonumber\\[.3pc]
b_{12} &= \in (C^{3}x^{3}+\dots +C^{r}x^{r}),\quad \in =\pm 1\nonumber\\[.3pc]
b_{iq} &= b_{pq}=0,\quad
\mathfrak{H}=\mathfrak{H}(x^{1},x^{2}),\dots,
\mathfrak{J}=\mathfrak{J}(x^{1},x^{2}).
\end{align}

Because of (2.1) and (2.2), the components of the first and the
second fundamental tensors of associate Bonnet hypersurface
$M^{\prime}$ can be written in the form:
\begin{align}
g_{11}^{\prime} &= g_{22}^{\prime}=(C^{3}x^{3}+\dots
+C^{r}x^{r})^{2}\frac{1}{\vert \mathfrak{J}\vert \sin
\theta}, \quad g_{pp}^{\prime}=1\nonumber\\[.3pc]
b_{11}^{\prime} &= (C^{3}x^{3}+\dots +C^{r}x^{r})
\frac{(\mathfrak{H}+ \mathfrak{J}\sin \theta)}{\vert
\mathfrak{J}\vert \sin \theta},\nonumber\\[.3pc]
b_{22}^{\prime} &= (C^{3}x^{3}+\dots
+C^{r}x^{r})\frac{(\mathfrak{H}- \mathfrak{J}\sin \theta)}{\vert
\mathfrak{J}\vert \sin\theta},\nonumber\\[.3pc]
b_{12}^{\prime} &= -\in (C^{3}x^{3}+\dots +C^{r}x^{r}),\quad
\in =\pm 1\nonumber\\[.3pc]
b_{iq}^{\prime} &= b_{pq}^{\prime}=0.
\end{align}

It is seen that all the fundamental quantities of $M$ and
$M^{\prime}$ are the same, except $b_{12}$ and $b_{12}^{\prime}$.
There is the relation $b_{12}=-b_{12}^{\prime}$ between $b_{12}$
and $b_{12}^{\prime}$.}
\end{case}

\begin{case}{\rm
Let us consider the case where all $T_{p}=0$. Then, it is easily
seen that in order to obtain the quantities in this case, it is
sufficient to replace $C^{3}x^{3}+\dots +C^{r}x^{r}$ by $1$ and
$(C^{3})^{2}+\dots +(C^{r})^{2}$ by $0$ in the quantities of
Case~$1$.

Thus, in this case, the fundamental quantities of $M$ and
$M^{\prime}$ in an orthogonal coordinate system
$x^{1},x^{2},\dots,x^{n}$ are respectively given by
\begin{align}
g_{11} &= g_{22}=\frac{1}{\vert \mathfrak{J}\vert \sin
\theta},\quad g_{pp}=1\nonumber\\[.3pc]
b_{11} &= \frac{(\mathfrak{H}+\mathfrak{J}\sin \theta)}{\vert \mathfrak{J}\vert \sin \theta},\nonumber\\[.3pc]
b_{22} &= \frac{(\mathfrak{H}-\mathfrak{J}\sin \theta)
}{\vert \mathfrak{J}\vert \sin \theta},\nonumber\\[.3pc]
b_{12} &= \in,\quad \in =\pm 1\nonumber\\[.3pc]
b_{iq} &= b_{pq}=0,\quad
\mathfrak{H}=\mathfrak{H}(x^{1},x^{2}),\dots,\quad
\mathfrak{J}=\mathfrak{J}(x^{1},x^{2})
\end{align}
and
\begin{align}
g_{11}^{\prime} &= g_{22}^{\prime}=\frac{1}{\vert
\mathfrak{J}\vert \sin \theta},\quad g_{pp}^{\prime}=1\nonumber\\[.3pc]
b_{11}^{\prime} &= \frac{(\mathfrak{H}+\mathfrak{J}\sin
\theta)}{\vert \mathfrak{J}\vert \sin \theta},\nonumber\\[.3pc]
b_{22}^{\prime} &= \frac{(\mathfrak{H}-\mathfrak{J}\sin
\theta)}{\vert \mathfrak{J}\vert \sin \theta},\nonumber\\[.3pc]
b_{12}^{\prime} &= -\in,\quad \in =\pm 1\nonumber\\[.3pc]
b_{iq}^{\prime} &= b_{pq}^{\prime}=0.
\end{align}}
\end{case}

\begin{defini}$\left.\right.$\vspace{.5pc}

\noindent{\rm If the fundamental quantities in an orthogonal
coordinate system $x^{1},x^{2},\dots,x^{n}$ are in the form
\begin{align*}
g_{11}=g_{22},\quad g_{pp}=1,\quad b_{12}=\in (C^{3}x^{3}+\dots
+C^{r}x^{r}),\quad b_{iq}=b_{pq}=0
\end{align*}
or
\begin{equation*}
g_{11}=g_{22},\quad g_{pp}=1,\quad b_{12}=\in,\quad
b_{iq}=b_{pq}=0,
\end{equation*}
where $C^{p}$ are constants. The net which consists of this
orthogonal coordinate system is called $A$-net.}
\end{defini}

Therefore, we have the following theorem.

\begin{theor}[({\it The Fundamental Theorem})] A~non-minimal hypersurface
in $R^{n+1}$ with no umbilical points is a Bonnet hypersurface if
and only if it has an $A$-net. Moreover{\rm ,} the fundamental
quantities of the Bonnet hypersurface and its Bonnet associate are
given by $(4.21)$ and $(4.22)$ respectively or $(4.23)$ and
$(4.24)$ respectively.
\end{theor}

\section*{Acknowledgement}

The authors would like to thank Markus Roussos for his interest in
this work.


\begin{thebibliography}{99}
\bibitem{1} Bobenko~A~I and Either~U, Bonnet surfaces and
Painlev\'e equations, {\it J~Reine Angew. Math.} {\bf 499} (1998)
47--79

\bibitem{2} Bonnet~O, Memoire sur la theorie des surfaces, Applicable
sur une surface donne, {\it J.~Ec. Poyt.} (1867) T. 25

\bibitem{3} Cartan~E, Sur les couples surfaces applicables avec
conservation des courbures principales, {\it Bull. Sci. Math.} {\bf 66}
(1942) 55--163

\bibitem{4} Chern~S~S, Deformations of surfaces preserving principal
curvature, Differantial Geometry and Complex Analysis, H.~E.~Rauch
Memorial Volume (eds) I~Chavel and H~M~Farkas (1985)
(Springer-Verlag) pp.~155--163

\bibitem{5} Chern~S~S, Bryant~R~L, Gardner~R~B, Goldschmidt~H~L and
Griffiths~P~A, Exterior Differential Systems, Mathematical Sciences
Research Institute Publications (1991) (New~York: Springer-Verlag)

\bibitem{6} Colares~A~G and Kenmotsu~K, Isometric deformation of
surfaces in preserving the mean curvature function, {\it
Pacific~J.~Math.} {\bf 136} (1989) 71--80

\bibitem{7} Csikos~B, Differential Geometry, Lectures Notes, Budapest
Semesters in Mathematics (1998)

\bibitem{8} Kenmotsu~K, An intrinsic characterization of $H$-deformable
surfaces, {\it J.~London Math. Soc.} {\bf 49} (1994) 555--568

\bibitem{9} Kobayashi~S and Nomizu~K, Foundations of differential
geometry (1969) (New~York: Interscience) vol.~2

\bibitem{10} Kokubu~M, Isometric deformations of hypersurfaces in a
{E}uclidean space preserving mean curvature, {\it Tohoku Math.~J.}
{\bf 44} (1992) 433--442

\bibitem{11} Roussos~I~M, Principal curvature preserving isometries of
surfaces in ordinary spaces, {\it Bol. Soc. Math.} {\bf 18(2)} (1987)
95--105

\bibitem{12} Roussos~I~M, Global results on Bonnet surfaces, {\it
J.~Geom.} {\bf 65} (1999) 151--158

\bibitem{13} Soyu\c{c}ok~Z, The problem of non-trivial isometries of
surfaces preserving principal curvatures, {\it J.~Geom.} {\bf 52} (1995)
173--188

\bibitem{14} Soyu\c{c}ok~Z, The problem of isometric deformations of a
Euclidean hypersurface preserving mean curvature, {\it Bull. Tech.
Univ.} {\bf 49} (1996) 551--562

\bibitem{15} Voss~K, Bonnet surfaces in spaces of constant curvature,
Lecture Notes~II of 1st MSJ Research Inst\i tute, Sendai, Japan (1993)
pp.~295--307

\bibitem{16} Xiuxiong~C and Chia-Kuei~P, Deformation of surfaces
preserving principal curvatures, {\it Lect. Notes Math.} (1989) 63--70
\end{thebibliography}
\end{document}